\documentclass[11pt]{article}
\usepackage{latexsym}
\usepackage{amsfonts}

\title{The spanning number and the independence number of a subset of an abelian group \\[.4in]}

\author{B\'{e}la Bajnok \\[.1in] Department of Mathematics, Gettysburg College \\
Gettysburg, PA 17325-1486 USA \\E-mail:  bbajnok@gettysburg.edu }

\date{April 29, 2003}

\newtheorem{thm}{Theorem}
\newtheorem{defin}[thm]{Definition}

\newtheorem{cor}[thm]{Corollary}
\newtheorem{prop}[thm]{Proposition}

\newcommand{\bc}[2]{{{#1}\choose{#2}}}

\begin{document}

\maketitle

\thispagestyle{empty}

\begin{abstract}

Let $A=\{a_1,a_2,\dots, a_m\}$ be a subset of a finite abelian group $G$.  We call $A$ {\it
$t$-independent} in $G$, if whenever $$\lambda_1a_1+\lambda_2a_2+\cdots +\lambda_m a_m=0$$ for some integers $\lambda_1, \lambda_2, \dots , \lambda_m$ with $$|\lambda_1|+|\lambda_2|+\cdots +|\lambda_m| \leq t,$$ we have $\lambda_1=\lambda_2= \cdots = \lambda_m=0$, and we say that $A$ is {\it $s$-spanning} in $G$, if every element $g$ of $G$ can be written as $$g=\lambda_1a_1+\lambda_2a_2+\cdots +\lambda_m a_m$$ for some integers $\lambda_1, \lambda_2, \dots , \lambda_m$ with $$|\lambda_1|+|\lambda_2|+\cdots +|\lambda_m| \leq s.$$

In this paper we give an upper bound for the size of a $t$-independent set and a lower bound for the size of an $s$-spanning set in $G$, and determine some cases when this extremal size occurs.  We also discuss an interesting connection to spherical combinatorics.

\end{abstract}

\section{Introduction}

We illuminate our concepts by the following examples.

{\bf Example 1}  Consider the set $A=\{1,4,6,9,11\}$ in the cyclic group $G=\mathbb{Z}_{25}$.  We are interested in the degree to which this set is \emph{independent} in $G$.  We find, for example, that $1+4+4-9=0$ and $11+11+9-6=0$, but that such an equation with only three terms from $A$ cannot be found.  We therefore say that $A$ is \emph{3-independent} in $G$ and write $\mathrm{ind}(A)=3$.  It can be shown that $A$ is optimal in each of the following regards:
\begin{itemize} 
\item no subset of $G$ of size $m>5$ is 3-independent in $G$ (furthermore, $A$ is essentially the unique 3-independent set in $G$ of size 5);
\item no subset of $G$ of size 5 is $t$-independent for $t>3$ (that is, for $t>3$, there will always be $t$, not necessarily distinct, elements with a signed sum of 0); and
\item $n=25$ is the smallest odd number for which a 3-independent set of size 5 in $\mathbb{Z}_{n}$ exists. (In fact, it can be shown that $\mathbb{Z}_{n}$ has a 3-independent set of size 5, if and only if, $n=20, 22, 24, 25, 26$, or $n \geq 28$.)
\end{itemize} 

The fact that $G=\mathbb{Z}_{25}$ has this relatively large 3-independent subset is due, as explained later, to the fact that $25$ has a prime divisor which is congruent to 5 mod 6.

{\bf Example 2}  How can one place a finite number of points on the $d$-dimensional sphere $S^d \subset \mathbb{R}^{d+1}$ with the highest \emph{momentum balance}?  For the circle $S^1$, the answer is given by the vertices of a regular polygon, but the issue is far more difficult for $d>1$.  For a positive integer $n$ and a set of integers $A=\{a_1,a_2,\dots, a_m\}$, define the set of $n$ points $X(A)=\{x_1,x_2,\dots, x_n\}$  with $$x_i=\frac{1}{\sqrt{m}} \cdot \left(\cos(\frac{2 \pi i a_1}{n}), \sin (\frac{2 \pi i a_1}{n}),  \dots, \cos (\frac{2 \pi i a_m}{n}), \sin (\frac{2 \pi i a_m}{n}) \right)$$ ( $i=1,2,\dots,n$); thus, for example, for $n=25$ and $A=\{1,4,6,9,11\}$, $X(A)$ is a set of 25 points on the unit sphere $S^9$.  It can be shown that this $X(A)$ is a \emph{spherical 3-design}, that is, for every polynomial $f: S^9 \rightarrow \mathbb{R}$ of total degree at most 3, the average value of $f$ on $S^9$ equals the arithmetic average of $f$ on $X(A)$.  We can also verify that $X(A)$ is optimal in that 
\begin{itemize}
\item no set of 25 points is a $t$-design on $S^9$ for $t>3$;
\item no set of 25 points is a 3-design on $S^d$ for $d >9$; 
\item $n=25$ is the minimum odd size for which a 3-design on $S^9$ exists.  (It was recently proved that an $n$-point 3-design on $S^9$ exists, if and only if, $n=20, 22, 24$, or $n \geq 25$.)
\end{itemize}

{\bf Example 3}  Finally, consider $A=\{3,4\}$ in $G=\mathbb{Z}_{25}$.  Note that every element of $G$ can be generated by a signed sum of at most three terms of $A$: $1=4-3, 2=3+3-4, \dots, 24=3-4$.  We therefore call $A=\{3,4\}$ a \emph{3-spanning set} in $G=\mathbb{Z}_{25}$, and write $\mathrm{span}(A)=3$.  Again, our example is extremal; it can be shown that 
\begin{itemize} 
\item no subset of $G$ of size $m<2$ is 3-spanning in $G$;
\item no subset of $G$ of size 2 is $s$-spanning for $s<3$; and
\item $n=25$ is the largest number for which a 3-spanning set of size 2 in $\mathbb{Z}_{n}$ exists. (Furthermore, as we will see, $\mathbb{Z}_{n}$ has a 3-spanning set of size at most 2 for every $n \leq 25$.)
\end{itemize} 

In fact, this example has an even more distinguished property: every element of $G$ can be written \emph{uniquely} as a signed sum of at most 3 elements of $A$; we call such a set \emph{perfect}.  As a consequence of being a perfect 3-spanning set, $A$ is also a maximum size 6-independent set in $G$.  The fact that $G=\mathbb{Z}_{25}$ has a perfect spanning subset of size 2 is due to the fact that 25 is the sum of two consecutive squares, as explained later. 

In the subsequent sections of this paper we define and investigate the afore-mentioned concepts and statements.  Topics similar to spanning numbers (e.g. $h$-bases) and independence numbers (e.g. sum-free sets, Sidon sets, and $B_h$ sequences) have been studied vigorously for a long time, see, for example, \cite{ErdFre:1991a}, \cite{Gra:1996a}, \cite{HalRot:1983a}, \cite{Nat:1996a}, \cite{Ruz:1993a}, \cite{Ruz:1995a}, \cite{WalStrWal:1972a}, and various sections of Guy's book \cite{Guy:1994a}.   For general references on spherical designs, see \cite{Ban:1988a}, \cite{DelGoeSei:1977a}, \cite{God:1993a}, \cite{GoeSei:1979a}, \cite{GoeSei:1981a}, \cite{Hog:1996a}, \cite{Rez:1992a}, and \cite{Sei:1996a}.

\section{Spanning numbers}

Let $G$ be a finite abelian group of order $|G|=n$, written in additive notation.  We are interested in the degree to which a given subset of $G$ spans $G$.  More precisely, we introduce the following definition.

\begin{defin} \label{s-span}

Let $s$ be a non-negative integer and $A=\{a_1,a_2,\dots, a_m\}$.  We say that $A$ is an \textup{\textbf{$s$-spanning set}} in $G$, if every $g \in G$ can be written as $$g=\lambda_1a_1+\lambda_2a_2+\cdots +\lambda_m a_m$$ for some integers $\lambda_1, \lambda_2, \dots , \lambda_m$ with $$|\lambda_1|+|\lambda_2|+\cdots +|\lambda_m| \leq s.$$  We call the smallest $s$ for which $A$ is $s$-spanning the \textup{\textbf{spanning number}} of $A$ in $G$, and denote it by $\mathrm{span}(A)$.

\end{defin}

Equivalently, $A$ is an $s$-spanning subset of $G$ if for every element $g \in G$, we can find non-negative integers $h$ and $k$ and elements $x$ and $y$ in $G$, so that $x$ is the sum of $h$ (not necessarily distinct) elements of $A$, $y$ is the sum of $k$ (not necessarily distinct) elements of $A$, $h+k \leq s$, and $g=x-y$.

The case $s=0$ is trivial: the only group $G$ which has a 0-span is the one with a single element; therefore, we may assume that $s \geq 1$ and $n \geq 2$.  Obviously, $A=G$ is an $s$-spanning subset of $G$ for every $s \geq 1$.  Here we are interested in small $s$-spanning sets in $G$; we denote the size of a minimum $s$-spanning set of $G$ by $p(G,s)$.

For $s=1$, it is clear that $\mathrm{span}(A)=1$, if and only if, for each $g \in G$, $A$ contains at least one of $g$ or $-g$; in particular, $A$ must contain every element of order 2.  Let $O(G,2)$ denote the set of order 2 elements of $G$; with this notation we have 

\begin{equation} \label{p(G,1)}  
p(G,1)=|O(G,2)|+ \frac{|G \setminus O(G,2) \setminus \{0\}|}{2}= \frac{n+|O(G,2)|-1}{2}.  
\end{equation}
As a special case, for the cyclic group of order $n$ we have 
\begin{equation} \label{p(Z,1)}  
p(\mathbb{Z}_n,1)=\lfloor n/2 \rfloor.  
\end{equation}

For $s \geq 2$, values of $p(G,s)$ seem difficult to establish, even in the case of the cyclic groups.  Computational data shows that  
\begin{equation} \label{p2}
p(\mathbb{Z}_n, 2) = \left\{
\begin{array}{cl}
0 & \mbox{if $n={\bf 1};$}\\
1 & \mbox{if $n=2, 3, 4, {\bf 5};$}\\
2 & \mbox{if $n=6, 7, \dots, 12, {\bf 13};$}\\
3 & \mbox{if $n=14, 15, \dots, 21;$}\\
4 & \mbox{if $n=22, 23, \dots, 33$, and $n=35;$}\\
5 & \mbox{if $n=34$, $n=36, 37, \dots, 49$, and $n=51;$}\\
\end{array}\right.
\end{equation}
and 
\begin{equation} \label{p3}
p(\mathbb{Z}_n, 3) = \left\{
\begin{array}{cl}
0 & \mbox{if $n={\bf 1};$}\\
1 & \mbox{if $n=2, 3, \dots, 6, {\bf 7};$}\\
2 & \mbox{if $n=8, 9 \dots, 24, {\bf 25};$}\\
3 & \mbox{if $n=26, 27, \dots, 50$, $n=52$, and $n=55;$}\\
4 & \mbox{if $n=51,53, 54$, $n=56, 57, \dots, 100$, and $n=104.$}\\
\end{array}\right.
\end{equation}
(Values marked in bold-face will be discussed later.)

As these values indicate, $p(\mathbb{Z}_n, s)$ is, in general, not a monotone function of $n$, though we believe that $$P(s):= \lim_{n \rightarrow \infty} \frac{p(\mathbb{Z}_n,s)^{s}}{n}$$ exists for every $s$.  The following theorem provides a lower bound for $p(G,s)$ which is of the order $n^{1/s}$ as $n$ goes to infinity. 

\begin{thm} \label{span-bound}

Let $m$ and $s$ be positive integers, and define $a(m,s)$ recursively by $a(m,0)=a(0,s)=1$ and $$a(m,s)=a(m-1,s)+a(m,s-1)+a(m-1,s-1).$$  

\begin{enumerate}

\item We have $$a(m,s)=\sum_{k=0}^s {s \choose k} {m \choose k} 2^k.$$

\item If $G$ has order $n$ and contains an $s$-spanning set of size $m$, then $n \leq a(m,s)$.

\end{enumerate}

\end{thm}

{\it Proof.}  1.  Let us define $$a'(m,s):=\sum_{k=0}^s {s \choose k} {m \choose k} 2^k.$$  Clearly, $a'(m,0)=a'(0,s)=1$; below we prove that $a'(m,s)$ also satisfies the recursion.

We have 
\begin{eqnarray*}
a'(m-1,s-1) & = & \sum_{k=0}^{s-1} {s-1 \choose k} {m-1 \choose k} 2^k \\
& = & \sum_{k=0}^{s-2} {s-1 \choose k} {m-1 \choose k} 2^k +{m-1 \choose s-1} 2^{s-1}, 
\end{eqnarray*}
and
\begin{eqnarray*}
a'(m-1,s) & = & \sum_{k=0}^{s} {s \choose k} {m-1 \choose k} 2^k \\
& = & \sum_{k=0}^{s-1} {s \choose k} {m-1 \choose k} 2^k +{m-1 \choose s} 2^s \\ 
& = &  \sum_{k=0}^{s-1} {s-1 \choose k-1} {m-1 \choose k} 2^k+\sum_{k=0}^{s-2} {s-1 \choose k} {m-1 \choose k} 2^k +{m-1 \choose s-1} 2^{s-1} +{m-1 \choose s} 2^s. 
\end{eqnarray*}

Next, we add $a'(m-1,s)$ and $a'(m-1,s-1)$.  Note that 
$${m-1 \choose s-1} 2^{s-1} +{m-1 \choose s-1} 2^{s-1}+{m-1 \choose s} 2^s={m \choose s} 2^{s},$$
and $$\sum_{k=0}^{s-2} {s-1 \choose k} {m-1 \choose k} 2^k + \sum_{k=0}^{s-2} {s-1 \choose k} {m-1 \choose k} 2^k=\sum_{k=0}^{s-2} {s-1 \choose k} {m-1 \choose k} 2^{k+1},$$ and by replacing $k$ by $k-1$, this sum becomes 
$$\sum_{k=0}^{s-1} {s-1 \choose k-1} {m-1 \choose k-1} 2^k.$$

Therefore, 
\begin{eqnarray*}
a'(m-1,s)+a'(m-1,s-1) & = & \sum_{k=0}^{s-1} {s-1 \choose k-1} {m-1 \choose k} 2^k+\sum_{k=0}^{s-1} {s-1 \choose k-1} {m-1 \choose k-1} 2^k+{m \choose s} 2^{s} \\
& = & \sum_{k=0}^{s-1} {s-1 \choose k-1} {m \choose k} 2^k+{m \choose s} 2^{s} \\
& = & \sum_{k=0}^{s} {s-1 \choose k-1} {m \choose k} 2^k \\ 
& = & \sum_{k=0}^{s} {s \choose k} {m \choose k} 2^k - \sum_{k=0}^{s} {s-1 \choose k} {m \choose k} 2^k \\
& = & a'(m,s)-a'(m,s-1).
\end{eqnarray*}

2.  Assume that $A=\{a_1,\dots,a_m \}$ is an $s$-spanning set in $G$ of size $m$, and let $$\Sigma= \{\lambda_1a_1+\cdots +\lambda_m a_m \mbox{    } | \mbox{    } \lambda_1, \dots , \lambda_m \in \mathbb{Z}, |\lambda_1|+\cdots +|\lambda_m| \leq s \}.  $$
We will count the elements in the index set $$I=\{(\lambda_1,\cdots ,\lambda_m) \mbox{    } | \mbox{    } \lambda_1, \dots , \lambda_m \in \mathbb{Z}, |\lambda_1|+\cdots +|\lambda_m| \leq s \},$$ as follows.  For $k=0,1,2,\dots,m$, let $I_k$ be the set of those elements of $I$ where exactly $k$ of the $m$ co\"{o}rdinates are non-zero.  How many elements are in $I_k$?  We can choose which $k$ of the $m$ co\"{o}rdinates are non-zero in ${m \choose k}$ ways; w.l.o.g. let these co\"{o}rdinates be $\lambda_1, \lambda_2, \dots, \lambda_k$.  Next, we choose the values of $|\lambda_1|, |\lambda_2|, \dots, |\lambda_k|$: since the sum of these $k$ positive integers is at most $s$, we have ${s \choose k}$ choices.  Finally, each of these co\"{o}rdinates can be positive or negative, and therefore $$|I_k|={s \choose k} {m \choose k} 2^k,$$ and $$|I|=\sum_{k=0}^m {s \choose k} {m \choose k} 2^k=\sum_{k=0}^s {s \choose k} {m \choose k} 2^k=a(m,s).$$  Since $A$ is $s$-spanning in $G$, we must have $n =|\Sigma| \leq |I|=a(m,s)$.  $\quad \Box$

Theorem \ref{span-bound} thus provides a lower bound for the size of an $s$-spanning set in $G$ which is of the order $n^{1/s}$ as $n$ goes to infinity.  

For exact values, we establish the following results.

\begin{prop} \label{monotone} Let $s \geq 1$ be an integer.
\begin{enumerate}
\item If $2 \leq n \leq 2s+1$, then the set $\{1\}$ is $s$-generating in $\mathbb{Z}_n$ and $p(\mathbb{Z}_n,s)=1$.

\item If $2s+2 \leq n \leq 2s^2+2s+1$, then the set $\{s,s+1\}$ is $s$-generating in $\mathbb{Z}_n$ and $p(\mathbb{Z}_n,s)=2$.

\item If $n \geq 2s^2+2s+2$, then $p(\mathbb{Z}_n,s) \geq 3$.

\end{enumerate}

\end{prop}

{\it Proof.}  1 is trivial.  To prove 2, let $$\Sigma=\{\lambda_1s+\lambda_2(s+1) \mbox{    } | \mbox{    }  \lambda_1,\lambda_2 \in \mathbb{Z}, |\lambda_1|+|\lambda_2| \leq s \}.$$  The elements of $\Sigma$ lie in the interval $[-(s^2+s),(s^2+s)]$ and, since the index set $$I=\{(\lambda_1,\lambda_2) \mbox{    }  | \mbox{    }  \lambda_1,\lambda_2 \in \mathbb{Z}, |\lambda_1|+|\lambda_2| \leq s \}$$ contains exactly $2s^2+2s+1$ elements, it suffices to prove that no integer in $[-(s^2+s),(s^2+s)]$ can be written as an element of $\Sigma$ in two different ways.  Indeed, it is an easy exercise to show that $$\lambda_1s+\lambda_2(s+1)=\lambda_1's+\lambda_2'(s+1) \in \Sigma$$ implies $\lambda_1=\lambda_1'$ and $\lambda_2=\lambda_2'$; therefore, the set $\{s,s+1\}$ is $s$-generating in $\mathbb{Z}_n$.  As the $s$-span of a single element can contain at most $2s+1$ elements, for values $n \geq 2s+2$ we must have $p(\mathbb{Z}_n,s)=2$.  Statement 3 follows from Theorem \ref{span-bound} by noting that $a(2,s)=2s^2+2s+1$.  $\quad \Box$

Let us now examine the extremal cases of Theorem \ref{span-bound}. 

\begin{defin} \label{perfect}

Suppose that $A$ is an $s$-spanning set of size $m$ in $G$ and that $a(m,s)$ is defined as in Theorem \ref{span-bound}.  If $|G|=n=a(m,s)$, then we say that $A$ is a \textup{\textbf{perfect $s$-spanning set}} in $G$.

\end{defin}

Cases where $\mathbb{Z}_{n}$ has a perfect $s$-spanning set for $s=2$ and $s=3$ are marked with bold-face in (\ref{p2}) and (\ref{p3}).  Trivially, the empty-set is a perfect $s$-spanning set in $\mathbb{Z}_{1}$ for every $s$.  With (\ref{p(Z,1)}) and Proposition \ref{monotone}, we can exhibit some other perfect spanning sets in the cyclic group.

\begin{prop} \label{perfex}

Let $m$, $n$, and $s$ be positive integers, and let $G=\mathbb{Z}_{n}$.

\begin{enumerate}

\item If $n=2m+1$, then the set $\{1,2,\dots,m\}$ is a perfect 1-spanning set in $G$.

\item If $n=2s+1$, then the set $\{1\}$ is a perfect s-spanning set in $G$.

\item If $n=2s^2+2s+1$, then the set $\{s,s+1\}$ is a perfect s-spanning set in $G$.

\end{enumerate}

\end{prop}

Note that the sets given in Proposition \ref{perfex} are not unique: any element of the set in 1 can be replaced by its negative; in 2, the set $\{a\}$ is perfect for every $a$ which is relatively prime to $n$; it is not difficult to show that another example in 3 is provided by $A=\{1,2s+1\}$ (however, the set $\{s,s+1\}$ in Proposition \ref{monotone} cannot be replaced by $\{1,2s+1\}$).  We could not find perfect spanning sets for $s \geq 2$ and $m \geq 3$.  It might be an interesting problem to find and classify all perfect spanning sets.

\section{Independence numbers}

As in the previous section, we let $G$ be a finite abelian group of order $|G|=n$, written in additive notation, and suppose that $A$ is a subset of $G$.  Here we are interested in the degree to which $A$ is independent in $G$.  More precisely, we introduce the following definition.

\begin{defin} \label{t-independent}

Let $t$ be a non-negative integer and $A=\{a_1,a_2,\dots, a_m\}$.  We say that $A$ is a \textup{\textbf{$t$-independent set}} in $G$, if whenever $$\lambda_1a_1+\lambda_2a_2+\cdots +\lambda_m a_m=0$$ for some integers $\lambda_1, \lambda_2, \dots , \lambda_m$ with $$|\lambda_1|+|\lambda_2|+\cdots +|\lambda_m| \leq t,$$ we have $\lambda_1=\lambda_2= \cdots = \lambda_m=0$.  We call the largest $t$ for which $A$ is $t$-independent the \textup{\textbf{independence number}} of $A$ in $G$, and denote it by $\mathrm{ind}(A)$.

\end{defin}

Equivalently, $A$ is a $t$-independent set in $G$, if for all non-negative integers $h$ and $k$ with $h+k \leq t$, the sum of $h$ (not necessarily distinct) elements of $A$ can only equal
the sum of $k$ (not necessarily distinct) elements of $A$ in a \emph{trivial} way, that is, $h=k$ and
the two sums contain the same terms in some order.

Here we are interested in the size of a maximum $t$-independent set in $G$; we denote this by $q(G,t)$.

Since $0 \leq \mathrm{ind}(A) \leq n-1$ holds for every subset $A$ of $G$ (so no subset is ``completely'' independent), we see that $q(G,0)=n$ and $q(G,n)=0$.  It is also clear that $\mathrm{ind}(A)=0$, if and only if, $0 \in A$, hence 
\begin{equation} \label{ob}  
q(G,1)=n-1.  
\end{equation}
For the rest of this section we assume that $t \geq 2$.  

We can easily determine the value of $q(G,2)$ as well.  First, note that $A$ cannot contain any element of $\{0\} \cup \mathrm{Ord}(G,2)$ (the elements of order at most 2); to get a maximum 2-independent set in $G$, take exactly one of each element or its negative in $G \setminus \mathrm{Ord}(G,2)\setminus \{0\}$, hence we have

\begin{equation} \label{obs}  
q(G,2)=\frac{n-|\mathrm{Ord}(G,2)|-1}{2}.  
\end{equation}
As a special case, for the cyclic group of order $n$ we have 
\begin{equation} \label{q(Z,2)}  
q(\mathbb{Z}_n,2)=\lfloor (n-1)/2 \rfloor.  
\end{equation}

Note that if $\mathrm{Ord}(G,2) \cup \{0\}=G$ then $q(G,2)=0$; for $n \geq 2$ this occurs only for the elementary abelian 2-group.  If $\mathrm{Ord}(G,2) \cup \{0\} \not = G$ then, since $\mathrm{Ord}(G,2)\cup \{0\}$ is a subgroup of $G$, we have $1 \leq |\mathrm{Ord}(G,2)|+1 \leq n/2$, and therefore we get the following.

\begin{prop}

If $G$ is isomorphic to the elementary abelian 2-group, then $q(G,2)=0$.  Otherwise $$\frac{1}{4}n \leq q(G,2) \leq \frac{1}{2}n.$$

\end{prop}

Let us now consider $t=3$.  As before, if $G$ does not contain elements of order at least 4, then $q(G,3)=0$; this occurs if and only if $G$ is isomorphic to the elementary abelian $p$-group for $p=2$ or $p=3$.  In \cite{BajRuz:2003a} we proved the following.

\begin{thm}[\cite{BajRuz:2003a}] \label{3bounds}

If $G$ is isomorphic to the elementary abelian $p$-group for $p=2$ or $p=3$, then $q(G,3)=0$.  Otherwise $$\frac{1}{9}n \leq q(G,3) \leq \frac{1}{4}n.$$

\end{thm}
These bounds can be attained since $q(\mathbb{Z}_9,3)=1$ and $q(\mathbb{Z}_4,3)=1$.  

For the cyclic group $\mathbb{Z}_n$, we can find explicit 3-independent sets as follows.  For every $n$, the odd integers which are less than $n/3$ form a 3-independent set; if $n$ is even, we can go up to (but not including) $n/2$ as then the sum of two odd integers cannot equal $n$.  We can do better in one special case when $n$ is odd; namely, when $n$ has a prime divisor $p$ which is congruent to 5 mod 6, one can show that the set 
\begin{equation} \label{5mod6}
\left\{ pi_1+2i_2+1 \mbox{    } | \mbox{    } i_1=0,1,\dots,\frac{n}{p}-1, \mbox{  } i_2=0,1,\dots,\frac{p-5}{6} \right\}
\end{equation}
is 3-independent.  It is surprising that these examples cannot be improved, as we have the following exact values.

\begin{thm}[\cite{BajRuz:2003a}] \label{3free}

For the cyclic group $G=\mathbb{Z}_n$ we have
$$q(\mathbb{Z}_n, 3) = \left\{
\begin{array}{cl}
\left\lfloor \frac{n}{4} \right\rfloor & \mbox{if $n$ is even,}\\
\left(1+\frac{1}{p}\right) \frac{n}{6} & \mbox{if $n$ is odd, has prime divisors congruent to 5 $\pmod 6$,} \\ & \mbox{and $p$ is the smallest such divisor,}\\
\left\lfloor \frac{n}{6} \right\rfloor & \mbox{otherwise.}\\
\end{array}\right.$$

\end{thm}

For $t \geq 4$, exact results seem more difficult.  
With the help of a computer, we generated the following values.

\begin{equation} \label{s4}
q(\mathbb{Z}_n, 4) = \left\{
\begin{array}{cl}
0 & \mbox{if $n={\bf 1}, 2, 3, 4;$}\\
1 & \mbox{if $n={\bf 5}, 6, \dots, 12;$}\\
2 & \mbox{if $n={\bf 13}, 14, \dots, 26;$}\\
3 & \mbox{if $n=27, 28, \dots, 45$, and $n=47;$}\\
4 & \mbox{if $n=46$, $n=48, 49, \dots, 68$, and $n=72, 73;$}\\
5 & \mbox{if $n=69, 70, 71$, and $n=74, 75, \dots, 102;$}\\
\end{array}\right.
\end{equation}
\begin{equation} \label{s5}
q(\mathbb{Z}_n, 5) = \left\{
\begin{array}{cl}
0 & \mbox{if $n={\bf 1}, 2, 3, 4, 5;$}\\
1 & \mbox{if $n={\bf 6}, 7, \dots, 17$, and $n=19, 20$;}\\
2 & \mbox{if $n={\bf 18}$, $n=21, 22, \dots, 37$, $n=39, 40, 41$, $n=43, 44, 45, 47$;}\\
3 & \mbox{if $n={\bf 38}, 42, 46$, $n=48, 49, \dots, 69$, $n=71, 72, 73, 75, 76, 77, 79, 81, 83, 85, 87;$}\\
\end{array}\right.
\end{equation}
and
\begin{equation} \label{s6}
q(\mathbb{Z}_n, 6) = \left\{
\begin{array}{cl}
0 & \mbox{if $n={\bf 1}, 2, 3, \dots, 6;$}\\
1 & \mbox{if $n={\bf 7}, 8, 9, \dots, 24$;}\\
2 & \mbox{if $n={\bf 25}, 26, 27, \dots, 69$;}\\
3 & \mbox{if $n=70, 71, \dots, 151$, and $n=153, 154, 155, 158, 159, 160.$}\\
\end{array}\right.
\end{equation}
(Values marked in bold-face will be discussed later.)

Again we see that $q(\mathbb{Z}_n, t)$ is not, in general, a monotone function of $n$; although for even values of $t$ the sequence seems to possess more regularity and we conjecture that $$Q(t):= \lim_{n \rightarrow \infty} \frac{q(\mathbb{Z}_n,t)^{t/2}}{n}$$ exists for every even $t$.  The following theorem establishes an upper bound for $q(G,s)$ which is of the order $n^{1/ \lfloor t/2 \rfloor}$ as $n$ goes to infinity. 

\begin{thm} \label{indep-bound}

Let $m$ and $t$ be positive integers, $t \geq 2$, and let us denote 
$$q(m,t) = \left\{
\begin{array}{cl}
a(m,t/2) & \mbox{if $t$ is even,}\\
a(m,(t-1)/2) +a(m-1,(t-1)/2)& \mbox{if $t$ is odd,}\\
\end{array}\right.$$
where $a(m,t)$ is defined in Theorem \ref{span-bound}.  If $G$ has order $n$ and contains a $t$-independent set of size $m$, then $n \geq q(m,t)$.

\end{thm}

{\it Proof.}  Assume that $A=\{a_1,\dots,a_m \}$ is a $t$-independent set in $G$ of size $m$, and define $$\Sigma= \{\lambda_1a_1+\cdots +\lambda_m a_m \mbox{    } | \mbox{    } \lambda_1, \dots , \lambda_m \in \mathbb{Z}, |\lambda_1|+\cdots +|\lambda_m| \leq \lfloor t/2 \rfloor \}$$ and $$I=\{(\lambda_1,\cdots ,\lambda_m) \mbox{    } | \mbox{    } \lambda_1, \dots , \lambda_m \in \mathbb{Z}, |\lambda_1|+\cdots +|\lambda_m| \leq \lfloor t/2 \rfloor \}.$$  

As in the proof of Theorem \ref{span-bound}, we have $|I|=a(m,\lfloor t/2 \rfloor)$.  Since $A$ is a $t$-independent set in $G$, the elements listed in $\Sigma$ must be all distinct, hence $n \geq |\Sigma| = |I|=a(m,\lfloor t/2 \rfloor)$.  If $t$ is even, we are done.

Now let $$\Sigma'= \{\lambda_1a_1+\cdots +\lambda_m a_m \mbox{    } | \mbox{    }  \lambda_1, \dots , \lambda_m \in \mathbb{Z}, \lambda_1 \geq 1, \lambda_1+|\lambda_2|+\cdots +|\lambda_m| = \lfloor t/2 \rfloor +1\}$$ and $$I'=\{(\lambda_1,\cdots ,\lambda_m) \mbox{    } | \mbox{    } \lambda_1, \dots , \lambda_m \in \mathbb{Z}, \lambda_1 \geq 1, \lambda_1+|\lambda_2|+\cdots +|\lambda_m| = \lfloor t/2 \rfloor +1 \}.$$  

We will count the elements in the index set $|I'|$ as follows.  For $k=0,1,2,\dots,m-1$, let $I_k$ be the set of those elements of $I'$ where exactly $k$ of the $m-1$ co\"{o}rdinates $\lambda_2, \dots , \lambda_m$ are non-zero.  An argument similar to that in the proof of Theorem \ref{span-bound} shows that $$|I_k'|={m-1 \choose k} {\lfloor t/2 \rfloor \choose k} 2^k,$$ hence $$|I'|=\sum_{k=0}^{m-1} {m-1 \choose k} {\lfloor t/2 \rfloor \choose k} 2^k=a(m-1,\lfloor t/2 \rfloor).$$
If $t$ is odd, then the elements listed in $\Sigma'$ must be distinct from each other and from those in $\Sigma$ as well, thus $n \geq |\Sigma|+|\Sigma'|=|I|+|I'|=a(m,\lfloor t/2 \rfloor)+a(m-1,\lfloor t/2 \rfloor).$  $\quad \Box$

Theorem \ref{indep-bound} thus provides an upper bound for the size of a $t$-independent set in $G$ which is of the order $n^{1/{\lfloor t/2 \rfloor}}$ as $n$ goes to infinity.  

For exact values, we establish the following results.

\begin{prop} \label{monot} Let $t \geq 2$ be an integer.
\begin{enumerate}
\item If $1 \leq n \leq t$, then $q(\mathbb{Z}_n,t)=0$.

\item If $t+1 \leq n \leq \lfloor t^2/2 \rfloor +t$, then the set $\{1\}$ is $t$-independent in $\mathbb{Z}_n$ and $q(\mathbb{Z}_n,t)=1$.

\item 

\begin{enumerate}

\item Suppose that $t$ is even.  If $n  \geq t^2/2+t+1$, then the set $\{t/2 , t/2 +1\}$ is $t$-independent in $\mathbb{Z}_n$ and $q(\mathbb{Z}_n,t) \geq 2$.

\item Suppose that $t$ is odd.  If $n  = (t^2-1)/2+t+1$, then the set $\{1,t\}$ is $t$-independent in $\mathbb{Z}_n$ and $q(\mathbb{Z}_n,t) = 2$.

\end{enumerate}

\end{enumerate}

\end{prop}

{\it Proof.}  Let $q(m,t)$ be defined as in Theorem \ref{indep-bound}. 
Since $q(1,t)=t+1$, our first claim follows from Theorem \ref{indep-bound}.  To prove 2, note that if $n \geq t+1$, then $\{1\}$ is $t$-independent in $\mathbb{Z}_n$; furthermore, $q(2,t)=\lfloor t^2/2 \rfloor +t+1$.

Now let $t$ be even, and assume that $n  \geq t^2/2+t+1$.  We define $$\Sigma=\{\lambda_1 \frac{t}{2}+\lambda_2(\frac{t}{2}+1) \mbox{    } | \mbox{    }  \lambda_1,\lambda_2 \in \mathbb{Z}, |\lambda_1|+|\lambda_2| \leq t \}.$$  The elements of $\Sigma$ lie in the interval $[-(t^2/2+t),(t^2/2+t)]$ and therefore, to prove 3 (a), it suffices to show that $$0 = \lambda_1 \frac{t}{2}+\lambda_2(\frac{t}{2}+1) \in \Sigma$$ implies $\lambda_1=\lambda_2=0$, which is an easy exercise.  Statement 3 (b) is essentially similar.
  $\quad \Box$

We now turn to the extremal cases of Theorem \ref{indep-bound}. 

\begin{defin} \label{tight}

Suppose that $A$ is a $t$-independent set of size $m$ in $G$ and that $q(m,t)$ is defined as in Theorem \ref{indep-bound}.  If $|G|=n=q(m,t)$, then we say that $A$ is a \textup{\textbf{tight $t$-independent set}} in $G$.

\end{defin}

Cases where $\mathbb{Z}_{n}$ has a tight $t$-independent set for $t=4$, $t=5$, and $t=6$ are marked with bold-face in (\ref{s4}), (\ref{s5}), and (\ref{s6}).  Trivially, the empty-set is a perfect $t$-independent set in $\mathbb{Z}_{1}$ for every $t$.  With (\ref{ob}), (\ref{q(Z,2)}), Theorem \ref{3free}, Proposition \ref{monot}, and one other (sporadic) example, we have the following tight $t$-independent sets in the cyclic group.

\begin{prop} \label{tightex}

Let $m$, $n$, and $t$ be positive integers, and let $G=\mathbb{Z}_{n}$.

\begin{enumerate}

\item If $n=2$, then the set $\{1\}$ is a tight 1-independent set in $G$.

\item If $n=2m+1$, then the set $\{1,2,\dots,m\}$ is a tight 2-independent set in $G$.

\item If $n=4m$, then the set $\{1,3,\dots,2m-1\}$ is a tight 3-independent set in $G$.

\item If $n=t+1$, then the set $\{1\}$ is a tight $t$-independent set in $G$.

\item Let $n=\lfloor t^2/2 \rfloor +t+1$.  If $t$ is even, then the set $\{t/2 ,t/2+1\}$ is a tight $t$-independent set in $G$; if $t$ is odd, then the set $\{1,t\}$ is a tight $t$-independent set in $G$.

\item If $n=38$, then the set $\{1, 7, 11\}$ is a tight $5$-independent set in $G$.

\end{enumerate}

\end{prop}

Proposition \ref{tightex} contains every tight (non-empty) $t$-independent set that we could find so far; in particular, we could not find tight $t$-independent sets for $t \geq 4$ and $m \geq 3$ other than the seemingly sporadic example listed last.  The problem of finding and classifying all tight $t$-independent sets remains open.

As it is clear from our exposition, there is a strong relationship between $s$-spanning sets and $t$-independent sets when $t$ is even.  Namely, we have the following.

\begin{thm}

Let $s$ and $t$ positive integers, $t$ even.  Let $A$ be a subset of $G$, and suppose that $\mathrm{span}(A)=s$ and $\mathrm{ind}(A)=t$.

\begin{enumerate}

\item The order $n$ of $G$ satisfies $a(m,t/2) \leq n \leq a(m,s)$.

\item We have $t \leq 2s$.

\item The following three statements are equivalent.

\begin{enumerate}
\item[(i)] $t=2s$;

\item[(ii)] $A$ is a perfect $s$-spanning set in $G$; and

\item[(iii)] $A$ is a tight $t$-independent set in $G$.

\end{enumerate}

\end{enumerate} 

\end{thm}

The analogous relationship when $t$ is odd is considerably more complicated and will be the subject of future study.

\section{Spherical designs}

Here we discuss an application of the previous section to spherical combinatorics.  We are interested in placing a finite number of points on the $d$-dimensional sphere $S^d \subset \mathbb{R}^{d+1}$ with the highest \emph{momentum balance}.   The following definition was introduced by Delsarte, Goethals, and Seidel in 1977 \cite{DelGoeSei:1977a}. 

\begin{defin} \label{designs}

Let $t$ be a non-negative integer.  A finite set $X$ of points on the $d$-sphere $S^d \subset \mathbb{R}^{d+1}$ is a \textup{\textbf{spherical $t$-design}}, if for every polynomial $f$ of total degree $t$ or less, the average value of $f$ over the whole sphere is equal to the arithmetic average of its values on $X$.

\end{defin}

In other words, $X$ is a spherical $t$-design if the Chebyshev-type quadrature formula
\begin{equation} \label{eq:quad}
\frac{1}{\sigma_d(S^d)} \int_{S^d}f({\bf x}) d \sigma_d({\bf x}) \approx \frac{1}{|X|} \sum_{{\bf x} \in X} f({\bf x})
\end{equation}
is exact for all polynomials $f: S^d \rightarrow \mathbb{R}$ of total degree at most $t$ ($\sigma_d$ denotes the surface measure on $S^d$).

The concept of $t$-designs on the sphere is analogous to $t-(v,k,\lambda)$ designs in combinatorics (see \cite{Sei:1990a}), and has been studied in various contexts, including representation theory, spherical geometry, and approximation theory.  For general references see \cite{Ban:1988a}, \cite{DelGoeSei:1977a}, \cite{God:1993a}, \cite{GoeSei:1979a}, \cite{GoeSei:1981a}, \cite{Hog:1996a}, \cite{Rez:1992a}, and \cite{Sei:1996a}.  The existence of spherical designs for every $t$ and $d$ and large enough $n=|X|$ was first proved by Seymour and Zaslavsky in 1984 \cite{SeyZas:1984a}.
 
A central question in the field is to find all integer triples $(t,d,n)$ for which a spherical $t$-design on $S^d$ exists consisting of $n$ points, and to provide explicit constructions for these parameters.  Clearly, to achieve high momentum balance on the sphere, one needs to take a large number of points.  Delsarte, Goethals, and Seidel \cite{DelGoeSei:1977a} provide the tight lower bound
\begin{equation} \label{eq:lower}
n \geq N_t^d:=\bc{d+\lfloor t/2 \rfloor}{\lfloor t/2 \rfloor}+\bc{d+\lfloor (t-1)/2 \rfloor}{\lfloor (t-1)/2}. 
\end{equation}
We shall refer to the bound $N_t^d$ in (\ref{eq:lower}) as the DGS bound.  Spherical designs of this minimum size are called \textup{\textbf{tight}}.  Bannai and Damerell \cite{BanDam:1979a}, \cite{BanDam:1980a} proved that tight spherical designs for $d \geq 2$ exist only for $t=1,2,3,4,5,7$, or $11$. All tight $t$-designs are known, except possibly for $t=4,5,$ or 7. In particular, there is a unique $11-$design ($d$=23 and $n=196560$).

Let us now attempt to construct spherical designs.  One's intuition that the vertices of a regular polygon provide spherical designs on the circle $S^1$ is indeed correct; more precisely, we have the following.

\begin{prop} \label{d=1}

Let $t$ and $n$ be positive integers.
\begin{enumerate}
\item
If $n \leq t$, then there is no $n$-point spherical $t$-design on $S^1$.
\item
Suppose that $n \geq t+1$.  For a positive integer $j$, define 
\begin{equation} \label{z_n^j}
{\bf z}_n^j:=\left(\cos(\frac{2 \pi j}{n}), \sin (\frac{2 \pi j }{n}) \right).
\end{equation}  Then the set $X_n:=\{{\bf z}_n^j | j=1,2,\dots,n\}$ is a $t$-design on $S^1$.
\end{enumerate}
\end{prop} 

{\it Proof.}  1 follows from the DGS bound as $N_t^1=t+1$.  To prove 2, we first note that, using spherical harmonics, one can prove (see \cite{DelGoeSei:1977a}) that, in general, a finite set $X$ is a spherical $t$-design, if and only if, for every integer $1 \leq k \leq t$ and every homogeneous \emph{harmonic} polynomial $f$ of total degree $k$, $$\sum_{{\bf x} \in X} f({\bf x})=0.$$  (A polynomial is harmonic if it is in the kernel of the Laplace operator.)  The set of homogeneous harmonic polynomials of total degree $k$ on the circle, $\mathrm{Harm}_k(S^1)$, is a 2-dimensional vector space over the reals and is spanned by the polynomials $\mathrm{Re}(z^k)$ and $\mathrm{Im}(z^k)$ where $z=x+\sqrt{-1}y$ (we can think of the elements of $X$ and $S^1$ as complex numbers).  Therefore, we see that $X$ is a $t$-design on $S^1$, if and only if, $$\sum_{{\bf z} \in X} {\bf z}^k=0$$ for $k=1,2,\dots,t$.  With $X_n$ as defined above, one finds that $$\sum_{j=1}^n ({\bf z}_n^j)^k =\left\{ \begin{array}{ll}
0 & \mbox{if $k \not \equiv 0$ mod $n$}, \\ \\
n & \mbox{if $k \equiv 0$ mod $n$}. \end{array}
\right.$$
Therefore, $X_n$ is a $t$-design on $S^1$, if and only if, $k \not \equiv 0$ mod $n$ for $k=1,2,\dots, t$ (using the terminology of our last section, if and only if, $\{1\}$ is a $t$-independent set in $\mathbb{Z}_{n}$), or $n\geq t+1$.  $\quad \Box$

A further classification of $t$-designs on the circle can be found in Hong's paper \cite{Hon:1982a}; he proved, for example, that if $n \geq 2t+3$, then there are infinitely many $t$-designs on $S^1$ which do not come from regular polygons.  

We now attempt to generalize Proposition \ref{d=1} for higher dimensions.  For simplicity, we assume that $d$ is odd, and let $m=(d+1)/2$.  (The case when $d$ is even can be reduced to this case by a simple technique, see \cite{Baj:1998a} or \cite{Mim:1990a}.)

Let $a_1, a_2, \dots, a_m$ be integers, and set $A:=\{a_1, a_2, \dots, a_m \}$.  For a positive integers $n$, define \begin{equation} \label{X_n(A)} X_n(A):=\left\{ \frac{1}{\sqrt{m}} \left( {\bf z}_n^j(a_1),  {\bf z}_n^j(a_2), \dots, {\bf z}_n^j(a_m) \right) \mbox{   } | \mbox{   }  j=1,2,\dots,n \right\},\end{equation}
where, like in (\ref{z_n^j}), $${\bf z}_n^j(a_i):=\left( \cos (\frac{2 \pi j}{n}a_i), \sin (\frac{2 \pi j}{n}a_i) \right).$$

Note that $X_n(A) \subset S^d$.  In \cite{Baj:1998a} we proved the following. 

\begin{thm}[\cite{Baj:1998a}] \label{connect}

Let $t$, $d$, and $n$ be positive integers with $t \leq 3$, $d$ odd, and set $m=(d+1)/2$.  For integers $a_1, a_2, \dots, a_m$, define $X_n(A)$ as in (\ref{X_n(A)}).  If $A$ is a $t$-independent set in $\mathbb{Z}_{n}$, then $X_n(A)$ is a spherical $t$-design on $S^d$.  

\end{thm}

Theorem \ref{connect} yields the following results.  

\begin{cor} \label{1,2,3}

Let $n$ and $d$ be positive integers, $d$ odd, and set $m=(d+1)/2$.

\begin{enumerate}

\item 
\begin{enumerate}

\item If $n=1$, then there is no $n$-point spherical $1$-design on $S^d$.
 
\item If $n \geq 2$, define $a_i=1$ for $1 \leq i \leq m$.  Then the set $X_n(A)$, as defined in (\ref{X_n(A)}), is a spherical $1$-design on $S^d$. 

\end{enumerate}

\item 
\begin{enumerate}

\item If $n \leq d+1$, then there is no $n$-point spherical $2$-design on $S^d$.
 
\item If $n \geq d+2$, define $a_i=i$ for $1 \leq i \leq m$.  Then the set $X_n(A)$, as defined in (\ref{X_n(A)}), is a spherical $2$-design on $S^d$. 

\end{enumerate}

\item 
\begin{enumerate}

\item If $n \leq 2d+1$, then there is no $n$-point spherical $3$-design on $S^d$.
 
\item If $n \geq 2d+2$ is even or if $n \geq 3d+3$ is odd, define $a_i=2i+1$ for $1 \leq i \leq m$; if $$n \geq \frac{p}{p+1}(3d+3)$$ where $p$ is a divisor of $n$ which is congruent to 5 mod 6, choose $A$  to be any $m$ elements of the set in (\ref{5mod6}).  In each case the set $X_n(A)$, as defined in (\ref{X_n(A)}), is a spherical $3$-design on $S^d$. 

\end{enumerate}
\end{enumerate}  

\end{cor}

{\it Proof.}  Parts (a) are from the DGS bounds $N_t^d$ for $t \leq 3$; parts (b) follow from Theorem \ref{connect} since, by (\ref{ob}), (\ref{q(Z,2)}), and the paragraph before Theorem \ref{3free}, the sets specified are $t$-independent for $t=1, 2$, and $3$, respectively (note that in all cases of 2 and 3, $m=(d+1)/2 \leq q(\mathbb{Z}_{n},t)$).  $\quad \Box$

Part 3 of Corollary \ref{1,2,3} leaves the question of existence of $3$-designs open for some odd values of $n$.  Note that the minimum value of $$\frac{p}{p+1}(3d+3)$$ is $5(d+1)/2$ (when $n$ is divisible by 5).  In \cite{Baj:1998a} we proved that a spherical 3-design on $S^d$ ($d$ odd) exists for \emph{every} odd value of $n \geq 5(d+1)/2$, and conjectured that 3-designs do not exist with $2(d+1) < n < 5(d+1)/2$ and $n$ odd.  This conjecture is supported by the numerical evidence of Hardin and Sloane \cite{HarSlo:1996a}.  A recent result of Boumova, Boyvalenkov, and Danev \cite{BouBoyDan:2002a} proves that no 3-design exists of odd size $n$ with $n< \approx 2.32(d+1)$.  In particular, the case $d=9$ of Example 2 in our Introduction is completely settled: 3-designs on $n$ points on $S^9$ exist, if and only if, $n \geq 20$ even, or $n \geq 25$ odd.
 
The application of $t$-independent sets to spherical $t$-designs seems more complicated when $t \geq 4$, and will be the subject of an upcoming paper.

{\bf Acknowledgments.} The author expresses his gratitude to his students Nicolae Laza for valuable computations and Nikolay Doskov for an improvement of Proposition \ref{monotone}.


\begin{thebibliography}{10}

\bibitem{Baj:1992a}
B.~Bajnok.
\newblock Construction of spherical $t$-designs.
\newblock {\em Geom. Dedicata}, 43:167--179, 1992.

\bibitem{Baj:1998a}
B. Bajnok.
\newblock Constructions of spherical 3-designs.
\newblock {\em Graphs Combin.}, 14/2:97--107, 1998.

\bibitem{BajRuz:2003a}
B. Bajnok and I. Ruzsa.
\newblock The independence number of a subset of an abelian group.
\newblock {\em Integers}, 3/Paper A2, 23 pp. (electronic), 2003.

\bibitem{Ban:1988a}
E.~Bannai.
\newblock On extremal finite sets in the sphere and other metric spaces.
\newblock {\em London Math. Soc. Lecture Note Ser.}, 131:13--38, 1988.

\bibitem{BanDam:1979a}
E.~Bannai and R.~M. Damerell.
\newblock Tight spherical designs {I}.
\newblock {\em J. Math. Soc. Japan}, 31:199--207, 1979.

\bibitem{BanDam:1980a}
E.~Bannai and R.~M. Damerell.
\newblock Tight spherical designs {II}.
\newblock {\em J. London Math. Soc. (2)}, 21:13--30, 1980.

\bibitem{BouBoyDan:2002a}
S. ~Boumova, P.~Boyvalenkov, and D.~Danev.
\newblock New nonexistence results for spherical designs.
\newblock In B. Bojanov, editor, {\em Constructive Theory of Functions}, pages 225--232, Varna, 2002.

\bibitem{DelGoeSei:1977a}
P.~Delsarte, J.~M. Goethals, and J.~J. Seidel.
\newblock Spherical codes and designs.
\newblock {\em Geom. Dedicata}, 6:363--388, 1977.

\bibitem{ErdFre:1991a}
P. Erd\H{o}s and R. Freud.
\newblock A {S}idon probl\'{e}mak{\"{o}}r.
\newblock {\em Mat. Lapok}, 1991/2:1--44, 1991.

\bibitem{God:1993a}
C.~D. Godsil.
\newblock {\em Algebraic Combinatorics}.
\newblock Chapman and Hall, Inc., 1993.

\bibitem{GoeSei:1979a}
J.~M. Goethals and J.~J. Seidel.
\newblock Spherical designs.
\newblock In D.~K. Ray-Chaudhuri, editor, {\em Relations between combinatorics
  and other parts of mathematics}, volume~34 of {\em Proc. Sympos. Pure Math.},
  pages 255--272. American {M}athematical {S}ociety, 1979.

\bibitem{GoeSei:1981a}
J.~M. Goethals and J.~J. Seidel.
\newblock Cubature formulae, polytopes and spherical designs.
\newblock In C.~Davis, B.~Gr{\"{u}}nbaum, and F.~A. Sher, editors, {\em The
  Geometric Vein: The Coveter Festschrift}, pages 203--218. Springer-{V}erlag
  New York, Inc., 1981.

\bibitem{Gra:1996a}
S.~W. Graham.
\newblock $B_h$ sequences.
\newblock In B.~C. Berndt, H.~G. Diamond, and A.~J. Hildebrand, editors, {\em Analytic number theory, Vol.1. (Allerton Park, IL, 1995)}, pages 431-449, Progr. Math. 138, Birkh\"{a}user Boston, Boston, MA, 1996.

\bibitem{Guy:1994a}
R.~K. Guy.
\newblock {\em Unsolved Problems in Number Theory}. Second edition.
\newblock Springer-Verlag New York, 1994.

\bibitem{HalRot:1983a}
H. Halberstam and K.~F. Roth.
\newblock {\em Sequences}. Second edition.
\newblock Springer-Verlag New York -- Berlin, 1983.

\bibitem{HarSlo:1996a}
R.~H. Hardin and N.~J.~A. Sloane.
\newblock Mc{L}aren's improved snub cube and other new spherical designs in
  three dimensions.
\newblock {\em Discrete Comput. Geom.}, 15:429--441, 1996.

\bibitem{Hog:1996a}
S.~G. Hoggar.
\newblock Spherical $t$-designs.
\newblock In C.~J. Colbourn and J.~H. Dinitz, editors, {\em The CRC handbook of
  combinatorial designs}, pages 462--466. CRC Press, Inc., 1996.

\bibitem{Hon:1982a}
Y. Hong.
\newblock On Spherical $t$-designs in $\mathbb{R}^2$.
\newblock {\em Europ. J. Combinatorics}, 3:255--258, 1982.

\bibitem{Mim:1990a}
J.~Mimura.
\newblock A construction of spherical 2-design.
\newblock {\em Graphs Combin.}, 6:369--372, 1990.

\bibitem{Nat:1996a}
M.~B. Nathanson.
\newblock {\em Additive Number Theory: Inverse Problems and the Geometry of Sumsets}.
\newblock Springer-Verlag New York, 1996.

\bibitem{Rez:1992a}
B.~Reznick.
\newblock Sums of even powers of real linear forms.
\newblock {\em Mem. Amer. Math. Soc.}, 463, 1992.

\bibitem{Ruz:1993a}
I. Ruzsa.
\newblock Solving linear equations in a set of integers I.
\newblock {\em Acta Arith.}, 65/3:259--282, 1993.

\bibitem{Ruz:1995a}
I. Ruzsa.
\newblock Solving linear equations in a set of integers II.
\newblock {\em Acta Arith.}, 72/4:385--397, 1995.

\bibitem{Sei:1990a}
J.~J. Seidel.
\newblock Designs and approximation.
\newblock {\em Contemp. Math.}, 111:179--186, 1990.

\bibitem{Sei:1996a}
J.~J. Seidel.
\newblock Spherical designs and tensors.
\newblock In E.~Bannai and A.~Munemasa, editors, {\em Progress in algebraic
  combinatorics}, volume~24 of {\em Adv. Stud. Pure Math.}, pages 309--321.
  Mathematical {S}ociety of {J}apan, 1996.

\bibitem{SeyZas:1984a}
P.~D. Seymour and T.~Zaslavsky.
\newblock Averaging sets: A generalization of mean values and spherical designs.
\newblock {\em Adv. Math.}, 52:213--240, 1984.

\bibitem{WalStrWal:1972a}
W.~D. Wallis, A.~P. Street, and J.~S. Wallis.
\newblock {\em Combinatorics: Room Squares, Sum-free Sets, Hadamard Matrices}, {\em Lecture Notes in Mathematics}, Vol. 292, Part 3.
\newblock Springer-Verlag, Berlin-New York, 1972.

\end{thebibliography}
\end{document}